\documentclass[journal]{IEEEtran}

\usepackage{cite}

\ifCLASSINFOpdf
   \usepackage[pdftex]{graphicx}
\else
\fi

\usepackage{amsmath}
\usepackage{xcolor}

\begin{document}

\title{Locational Marginal Value of Distributed\\ Energy Resources as Non-Wires Alternatives
}

\author{Panagiotis~Andrianesis,
        Michael~Caramanis,~\IEEEmembership{Senior~Member,~IEEE},
        Ralph~Masiello, ~\IEEEmembership{Life~Fellow,~IEEE},
        Richard~Tabors, ~\IEEEmembership{Member,~IEEE}, 
        and~Shay~Bahramirad, ~\IEEEmembership{Senior~Member,~IEEE}
        \thanks{The work of P. Andrianesis and M. Caramanis was supported in part by the Sloan Foundation under Grant G-2017-9723 and NSF AitF Grant 1733827. The work of R. Masiello and R. Tabors was supported in part by Commonwealth Edison (ComEd). P. Andrianesis and M. Caramanis are with the Systems Engineering Division, Boston University, Boston, MA (panosa@bu.edu, mcaraman@bu.edu). R. Masiello is with Quanta Technology, LLC., Raleigh, NC (rmasiello@quanta-technology.com). R. Tabors is with TCR Inc., Boston, MA (rtabors@tcr-us.com). S. Bahramirad is with ComEd, IL (Shay.Bahramirad@ComEd.com). }
}

\maketitle

\begin{abstract}
In this paper, we address the issue of valuating Distributed Energy Resources (DERs) as Non-Wires Alternatives (NWAs) against wires investments in the traditional distribution network planning process.
Motivated by the recent literature on Distribution Locational Marginal Prices, we propose a framework that allows the planner to identify rigorously the short-term Locational Marginal Value (LMV) of DERs using the notion of Marginal Cost of Capacity (MCC) of the best grid investment alternative to monetize hourly network constraint violations encountered during a yearly rate base timescale. 
We apply our methodology on two actual distribution feeders anticipated to experience overloads in the absence of additional DERs, and present numerical results on desirable LMV-based generic DER adoption targets and associated costs that can offset or delay different types of grid wires investments.
We close with a discussion on policy and actual DER adoption implementation.
\end{abstract}


\IEEEpeerreviewmaketitle

\section{Introduction}
\IEEEPARstart{D}{istribution} utilities have dealt with load growth by commensurate network investments. 
However, recent acceleration of Distributed Energy Resources (DERs) has raised the opportunity for considering DERs as Non-Wires Alternatives (NWAs) that enable deferral or avoidance of costly and often disruptive network investments.
In this vein, \cite{Tierney2016} eloquently posed a key question: what is the value of DERs to the distribution system? 
More specifically, what is the value of DERs at different hours and distribution network locations, and how do different DERs compare on an annual basis?

\subsection{Background and Motivation}

Traditionally, DERs referred to small and dispersed generation resources, such as solar or Combined Heat and Power (CHP), connected to the distribution network.
DERs were mainly associated with Distributed Generation (DG), whose value has been studied from various perspectives; \cite{WangNehrir2004, GeorgilakisHatziargyriou2013} consider optimal DG placement, \cite{LiZielke2005, GilJoos2006, WangEtAl2008, PiccoloSiano2009, WangEtAl2010, ZhangEtAl2016} the impact of DG on capacity deferral, and \cite{GroverSilvaEtAl2016} the grid's DER hosting capacity.
A sizable share of published work focuses on evaluating DG scenarios, their economics and impact on reliability and the environment.
Capacity deferral literature has so far relied among others on Avoided Cost, and Present Worth methods, (e.g., \cite{HeffnerEtAl1998, LiZielke2005}).
DG deferral of scheduled network upgrade investments for feeder groups are quantified in \cite{GilJoos2006}.
The reduction of power flowing over a radial feeder by adding DG is evaluated against the time it will take for the load to outgrow DG's effect.
The intuition that monetary benefits are maximized by adding DGs at the end of long feeders and near load pockets is confirmed in \cite{GilJoos2006}; 
the delay period till reinforcements are necessary is examined in \cite{WangEtAl2008}; \cite{PiccoloSiano2009} uses \cite{WangEtAl2008} to evaluate DG-related investment deferral value with feeder specific investments following \cite{GilJoos2006}.
Lastly, \cite{WangEtAl2010} quantifies the DG impact on demand growth and system security-related investments.

Although a widely acceptable definition of DERs is not yet cast in concrete, their concept has evolved to include not only DG (solar, CHP, small wind, etc.), but also energy storage, demand response, electric vehicles (EVs), microgrids, and energy efficiency.
Recently, estimating DER value by time and location is attracting increasing attention \cite{BodeEtAl2016, FineEtAl2016, RogersEtAl2017, HileEtAl2017, RobisonEtAl2017, DeMartiniEtAl2017, KristovEtAl2016}.
CA, NY, IL stakeholders have ongoing discussions on using the value of DERs as NWAs for compensation and incentive purposes.
In 2016, the California Public Utilities Commission approved a Locational Net Benefits Analysis framework \cite{CPUC2016}, and a Benefit Cost Analysis Framework was adopted in NY \cite{NY2016}.
Beyond the Brooklyn-Queens Demand Management demonstration project, NY utilities are announcing NWA projects and are actually procuring NWAs \cite{NYREV}.
The Illinois Commission of Commerce announced the Next Grid initiative in 2017 \cite{NextGrid} designating the value of DERs as a key focal point.

Despite the related literature, a consistent framework that compares DER adoption to traditional wires investments is still lacking. 
Indeed, in the current state-of-the-art, utility planners consider specific DERs assuming that their costs, capabilities, and the like, constitute known input to their NWA planning studies \cite{ContrerasEtAl2018, DeboeverEtAl2018}.
However, when the attraction of future DERs that are currently not in place is examined as a NWA, this input is in a state of flux, and hence unavailable with sufficient certainty.
Most importantly, since committing the study to uncertain input assumptions may affect its outcome significantly in favor or against specific technologies, regulators and stakeholders are likely, and justifiably so, to question them.
We propose a framework for considering DERs as NWAs that does not rely on \emph{guesses} of specific DER characteristics;
it is instead founded on quantifying \emph{generic} DER spatiotemporal marginal ``value-to-the-grid'' encompassing a marginal cost concept during hours of capacity constraint violations.

\subsection{Objectives and Contribution}

We strive to develop a formal framework that evaluates generic real and reactive power producing/consuming DERs as distribution NWAs.
High fidelity AC circuit analysis is used to estimate spatiotemporal marginal costs to the power system unbundled to their energy and grid components and quantify the generic DER spatiotemporal marginal value-to-the-grid.

The proposed framework builds upon short term locational marginal costing and pricing analysis \cite{SchweppeEtAl1988, CaramanisEtAl2016}. 
We rely upon and extend the concepts of the Marginal Cost of Capacity (MCC) and Locational Marginal Value (LMV) to quantify the value-to-the-grid of generic DER additions as NWAs that could or would be located on the grid to relieve constraint violations (e.g., line overloads, nodal over/under-voltages), while participating in available energy market products and services.
It should be noted that the terms LMV and MCC or similar expressions have been used in the literature of T{\&}D networks for several decades. 
For instance, Locational Marginal Prices (LMPs) characterize today's nodal electricity markets that originate from the seminal work on spot pricing of electricity \cite{SchweppeEtAl1988};
LMV has been used in a different context to characterize the value of storage capacity \cite{BoseBitar2014};
there is also an emerging literature on Distribution LMPs (DLMPs) \cite{CaramanisEtAl2016, YuanEtAl2018, BaiEtAl2018, Papavasiliou2018}.
The term of Marginal Distribution Capacity Cost (MDCC) has been also used extensively in the capacity deferral and DG planning literature \cite{HeffnerEtAl1998, WooEtAl1994, FengEtAl2009, Gutierrez2011}.
In this paper, LMV and MCC are construed differently to reflect the new context that they are used in.

More specifically, the MCC is computed from the cost of actual capital investments required to relieve anticipated constraint violations.
This cost is used to quantify the penalty for exacerbating constraints encountered in an infeasible AC OPF problem.
The LMV of a generic real power or reactive power DER represents the value of an incremental kW or kVAR provided to relieve the cost associated with violated constraints.
LMVs vary by node of the network and by hour.
As such, they assign values to specific DERs based on both their location and hourly profile across the year. 
Since our MCC computation results in a cost per unit of constraint violation, it impacts the LMV in a spatiotemporal manner to the extent that an incremental DER at a specific node and hour relieves each violated constraint with varying sensitivity.

The key contribution of our framework is that it 
\emph{(i)} relies on the cost of the best required wires investment to estimate generic kW and kVAR LMVs that are independent of any specific DER costs and capabilities, and
\emph{(ii)} provides the theoretically optimal amount and value of generic DERs required to defer the wires investment.
The associated annual DER procurement costs can be compared to the annual rate payer avoided costs that would have resulted from the deferred wires investment. Such comparisons performed on a yearly basis can inform whether DER adoption is a desirable non wires investment alternative.

Our approach embeds the explicit distribution planning problem into a spatiotemporal generic DER valuation framework, which is invariant of specific DER technologies and their associated costs.
Generic DER LMV is dependent only on the network characteristics, anticipated loads, constraint violations determined by detailed AC OPF, and the cost of required wires investments that may be needed to render the AC OPF problem feasible. 
Specific DERs required to alleviate network constraint violations can be construed as a composition of generic DER quantities.
The LMV of actual DERs and their affordable compensation can be derived from the generic DER LMV projected on actual DER potential real and reactive power hourly trajectories at their specific locations.

A high level introduction of the concepts of MCC and LMV in the context of the DER value-to-the-grid --- as applicable in our framework --- is presented in our preliminary work in \cite{TaborsEtAl2019}. 
In this paper, we thoroughly provide the network model and the method.
In addition, we apply our methodology on two representative test cases adapted from actual feeders of Commonwealth Edison (ComEd), IL, to examine DER-enabled deferral of wires investments.
The required wires investments include two typical although different cases of re-conductoring that were associated with anticipated load growth expected to result in line overload but no over/under-voltage violations.
Furthermore, we discuss several aspects of the resulting policy implications and extensions.

\subsection{Paper Organization}

The remainder of the paper is organized as follows. 
Section \ref{Model} presents the model formulation, Section \ref{Method} describes the framework from an algorithmic point of view, Section \ref{TestCases} introduces the test cases, and Section \ref{Results} presents numerical results.
Section \ref{Policy} discusses policy implications, and Section \ref{Conclusions} concludes and proposes future work.

\section{Model} \label{Model}

\begin{figure}[t]
\centering
\includegraphics[width=2.2in]{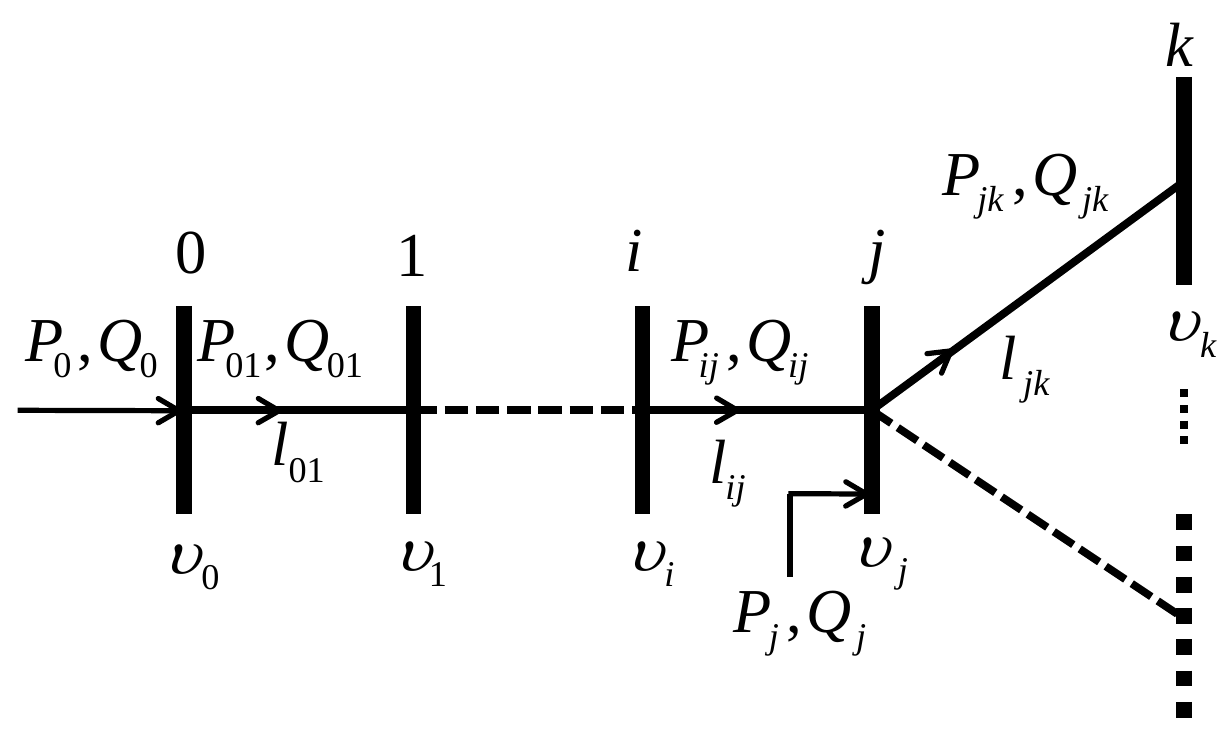}
\caption{Sketch of a tree network representation.}
\label{Fig1}
\end{figure}

We assume a balanced radial distribution network, represented by graph $(\mathcal{N},\mathcal{E})$. $\mathcal{N}$ is the set of nodes and $\mathcal{E}$ the set of edges.
Nodes are indexed by $0,1,...,n$, where $0$ is the root node. $\mathcal{N}\equiv \{0,1,...,n\}$, and ${\mathcal{N}^{+}}\equiv \mathcal{N}\backslash \{0\}$.
Pairs $(i,j)$ represent edges that denote lines connecting node $i$ with node $j$. 
The set of lines $\mathcal{E}$ has $n$ pairs, which are ordered by the $j$-th node.
The radial structure allows a unique path from the root node $0$ to node $j$, with $i$ the node that precedes $j$ in this path.
For each node $i\in \mathcal{N}$, let $V_i$ be the magnitude of the voltage, with $v_i \equiv V_{i}^2$, and minimum (maximum) voltage limits denoted by $V_{i}^{\min }$ ($V_{i}^{\max }$).
For each line $(i,j)\in E$, $r_{ij}$ is the resistance, $x_{ij}$ the reactance, $I_{ij}$ the magnitude of the current, with $l_{ij}\equiv I_{ij}^2$, $I_{ij}^{\max }$ the ampacity, and $P_{ij}$ and $Q_{ij}$ the sending-end real and reactive power flow, respectively.
$P_i$ and $Q_i$ denote the net real and reactive power injections at node $i$.
A positive (negative) value of $P_i$ refers to generation (consumption); similarly for the reactive power.
A sketch of a tree network is shown in Fig. \ref{Fig1}.

We use the DistFlow --- also referred to as the branch flow --- model introduced in \cite{BaranWu1989} and revised in \cite{FarivarLow2013}, which is a simplified, yet exact, representation of conventional AC power flow equations for a radial network.
The resulting AC OPF optimization problem is listed next.
\begin{equation}
\underset{P_0, Q_0, P_{ij}, Q_{ij}, v_i, l_{ij} } {\mathop{\min }} \,\,\,\,c^P P_0 + c^Q Q_0, \label{eq1}    
\end{equation}
subject to:
\begin{subequations}
\begin{align}
 P_0 = P_{01},\,\,\,\,\,\,\,(\lambda _{0}^{P}), \qquad \qquad \qquad \label{eq2a}\\
 Q_0 = Q_{01},\,\,\,\,\,\,\,(\lambda _{0}^{Q}),\qquad \qquad \qquad  \label{eq2b}\\
 P_{ij} - l_{ij} r_{ij} + P_j - \sum\limits_{k: j \to k} P_{jk} = 0, \,\,(\lambda _{j}^{P}) \,\,\,\, \forall j\in \mathcal{N}^{+}, \label{eq2c}	\\
 Q_{ij} - l_{ij} x_{ij} + Q_j - \sum\limits_{k: j \to k} Q_{jk} = 0, \,\,(\lambda _{j}^{Q}) \,\,\,\, \forall j\in \mathcal{N}^{+}, \label{eq2d}
\end{align}
\end{subequations}
\begin{equation}
v_j = v_i - 2(r_{ij} P_{ij} + x_{ij} Q_{ij} ) + ( r_{ij}^{2} + x_{ij}^{2} ) l_{ij}, \,\,\, \forall j\in \mathcal{N}^{+}, \label{eq3}     
\end{equation}
\begin{equation}
l_{ij} = \frac{ P_{ij}^{2} + Q_{ij}^{2} }{ v_i}, \quad \forall (i,j)\in \mathcal{E}, \label{eq4}  
\end{equation}
\begin{equation}
{{\left( V_{i}^{\text{min}} \right)}^{2}} \le v_i \le {{\left( V_{i}^{\text{max}} \right)}^{2}}, \quad \forall i\in \mathcal{N},  \label{eq5}  
\end{equation}
\begin{equation}
l_{ij} \le {{\left( I_{ij}^{\text{max}} \right)}^{2}}, \quad \forall (i,j)\in \mathcal{E},   \label{eq6}    
\end{equation}
where $P_0, Q_0, P_{ij}, Q_{ij}$ are real, and $v_i, l_{ij}$ non-negative. The time index is omitted for brevity.

The objective function \eqref{eq1} represents the cost of real and reactive power procured at the T\&D interface root node, with $c^P$ the real power LMP, and $c^Q$ a given reactive power compensation opportunity cost.
Notably, there is no transmission wholesale market price for reactive power, for reasons that, among others, include local market power concerns.
However, there is a cost for the provision of this service, which, in certain situations, can be viewed as the opportunity cost of a local generator (in the transmission system) providing this service, associated with foregoing the use of a unit of real power production \cite{CaramanisEtAl2016}.
We acknowledge that the issue of pricing/costing the provision of reactive power at the substation is complex, but further elaboration is beyond the scope of this paper.

The real and reactive power balance at each node are represented by \eqref{eq2a}--\eqref{eq2d}; their associated dual variables $\lambda _{i}^{P},\, \lambda _{i}^{Q}$ denote the real and reactive power DLMPs at node $i$. 
Constraints \eqref{eq3} and \eqref{eq4} define nodal voltage and line current. 
Constraints \eqref{eq5} and \eqref{eq6} impose voltage and current limits. We note that constraint \eqref{eq4} is non-convex. 
As proposed in \cite{FarivarLow2013}, replacing \eqref{eq4} by inequality
\begin{equation}
v_i l_{ij} \ge P_{ij}^{2} + Q_{ij}^{2},\quad \forall (i,j)\in \mathcal{E}, \label{eq7}   
\end{equation}
which is a convex Second Order Cone Programming (SOCP) constraint, introduces a convex relaxation of the problem.
For the cases of interest in this paper, this relaxation is exact; 
hence, instead of \eqref{eq4}, we will use \eqref{eq7} in our formulations.

\section{The Method} \label{Method}

The proposed framework's method consists of 3 steps: 
(1) Pre-processing (Subsection \ref{Preprocessing}), in which we calculate the constraint violation overload and the MCC; 
(2) Pricing (Subsection \ref{Pricing}), in which we obtain the real and reactive power LMVs for each hour and location, and 
(3) Generic DER Procurement (Subsection \ref{DERprocurement}), in which we derive the optimal addition of generic DERs that relieve the overload.

\subsection{Pre-processing} \label{Preprocessing}

\subsubsection{Overload Calculation}

Our aim is to calculate the amount of overload for each hour of the anticipated yearly load profile.
We employ the branch flow model, and we note that in the absence of inter-temporal constraints hourly calculations are parallelizable. 
In particular, omitting ampacity constraint \eqref{eq6}, we get the following OPF problem:
\begin{equation}
	\textbf{Opt1: } \eqref{eq1},\,\, \text{s.t. } \eqref{eq2a}-\eqref{eq2d}, \eqref{eq3}, \eqref{eq5}, \,\, \text{and } \eqref{eq7},  \label{eq8}  
\end{equation}
which, because of \eqref{eq7}, is a Quadratically Constrained Programming (QCP) problem, more specifically an SOCP problem.
We note that Opt1 essentially optimizes the voltage at the root node, since the net real/reactive power injections are fixed and the remaining variables (flows, currents, voltages) can be obtained by the load flow equations.
The solution of Opt1, which allows overload to occur, yields the values of $l_{ij,t}$, from which we calculate hourly overload $\Delta \hat{I}_{ij,t}$ in Amps for each line segment ($i,j$) exceeding its ampacity:
\begin{equation}
\Delta \hat{I}_{ij,t} = \max \left\{ 0,\sqrt{l_{ij,t}} - I_{ij}^{\text{max}} \right\}.	\label{eq9}    
\end{equation}
We use $\Delta \hat{I}_{ij,t}$ (and not $\Delta I_{ij,t}$) to distinguish the calculated (hat) values in the absence of the ampacity constraint \eqref{eq6}.

\subsubsection{MCC Calculation}
The MCC is calculated from the best grid investment cost, denoted by $C$ (in \$), obtained by a traditional wires solutions planning problem.

Let us first consider a case in which the best grid investment involves line upgrades, and hence the project cost $C$ can be directly allocated to each line segment.
Let $c_{ij}$ be the cost for increasing the line capacity (ampacity) by $\Delta I_{ij}^{\max }$ (in Amps), with $\sum\nolimits_{(i,j)}{c_{ij}} = C,$ and let $T_{ij}$ represent the number of hours in the year that the line is overloaded, i.e., the number of hours the line upgrade is required within the year. 
Since our horizon is one year, we first annualize the line upgrade cost to equal its anticipated impact on the rate base. 
For the purposes of this paper, we simply scale by a factor $\alpha$. We then define the MCC overload factor, denoted by $w_{ij}$, which we henceforth use interchangeably to MCC, as:
\begin{equation}
w_{ij} = \frac{\alpha \cdot c_{ij}}{ \Delta I_{ij}^{\max } \cdot T_{ij}}, \label{eq10}  
\end{equation}
where $w_{ij}$ (MCC) is measured in \$ per Amp of new capacity per (overloaded) hour, for the period of one year.\footnote{
This definition is in fact the average incremental cost of capacity. 
We use the term marginal for two reasons: (a) a small upgrade renders incremental an approximation of marginal, and (b) $w_{ij}$ is used in \eqref{eq12} as the coefficient of a linear ampacity overload cost where average and marginal coincide.}

Let us now consider a case in which the project involves an investment that cannot be allocated directly to the overloaded lines, e.g., building new lines as part of a reconfiguration scheme.
Arguably, we can still allocate the project cost to the overloaded lines, taking into account their maximum overload, $\Delta \hat{I}_{ij}^{\max } = \underset{t}{\mathop{\max }}\,\left\{ \Delta {{{\hat{I}}}_{ij,t}} \right\}$, and their length $L_{ij}$, as follows:
\begin{equation}
 c_{ij} = \frac{ \Delta \hat{I}_{ij}^{\max } L_{ij}} { \sum\nolimits_{(i,j)}{ \Delta \hat{I}_{ij}^{\max } L_{ij}}} C. \label{eq11}
\end{equation}
We can then apply \eqref{eq10} to derive the MCC, using the calculated value $\Delta \hat{I}_{ij}^{\max }$ instead of the actual increase in ampacity $\Delta I_{ij}^{\max }$ resulting from the line upgrade. 
Hence, we can view \eqref{eq11} as a reasonable, indirect, method for the allocation of the project cost, when a direct allocation is not applicable.

\subsection{Pricing} \label{Pricing}

In this step, we derive the generic DER spatiotemporal value. 
The idea is to monetize the overload $\Delta I_{ij,t}$ by the MCC factor $w_{ij}$; the new objective function that replaces \eqref{eq1} is:
\begin{equation}
\underset{P_0, Q_0, P_{ij}, Q_{ij}, v_i, l_{ij}, \Delta I_{ij}} {\mathop{\min }} \, c^P P_0 + c^Q Q_0 + \sum\limits_{(i,j)} w_{ij} \Delta I_{ij},  \label{eq12}    
\end{equation}
where the time index is omitted.
In \eqref{eq12}, $\Delta I_{ij}$ represents a new variable introduced for each overloaded line, so that the related costs are only applied to ($i,j$) exhibiting $\Delta I_{ij} > 0$ during a specific hour. 
Since the solution of Opt1 is known from the previous step, we define the overload variable $\Delta I_{ij}$ using the 1st order Taylor approximation, as follows:
\begin{equation}
\Delta I_{ij} = 0.5 { {\left( \sqrt{l_{ij}^{0}} \right)}^{-1}} l_{ij} + 0.5\sqrt{ l_{ij}^{0}} - I_{ij}^{\text{max}}, \label{eq13}
\end{equation}
where $l_{ij}^{0}$ is the current (magnitude squared) value derived from the solution of Opt1.

The cost for the overload in \eqref{eq12} represents the annualized pro-rated cost of the line, since we account only for the amount of new capacity needed in each hour, $\Delta I_{ij}$, instead of the maximum (lumpy) new capacity of the line $(\Delta I_{ij}^{\max })$.
Alternative approaches can be considered, as for instance, the Net Present Value of the annual revenue requirement of the capacity upgrade over an appropriate planning horizon. 
Our framework is applicable to such approaches, in fact, the subject of policy choices.
A key benefit is that the inclusion of the marginal avoided cost in $w_{ij}$ results in the DER investor and the customers sharing the avoided cost.
If the entire avoided cost of planned traditional investments, including excess capacity, were included in $w_{ij}$, then all of the avoided cost could be captured by generic DERs via the LMV mechanism, and customers/ratepayers would realize no net savings.

For each hour in which overload was identified in the solution of Opt1, we solve the following optimization problem:
\begin{equation}
	\textbf{Opt2: } \eqref{eq12},\,\, \text{s.t. } \eqref{eq2a}-\eqref{eq2d}, \eqref{eq3}, \eqref{eq5}, \eqref{eq7} \,\, \text{and } \eqref{eq13},  \label{eq14}  
\end{equation}
which is also a QCP (SOCP) problem. 
The LMVs are the shadow prices of \eqref{eq2c}--\eqref{eq2d}, i.e., $\lambda _{j}^{P}, \lambda _{j}^{Q}$, referred to as P-LMV and Q-LMV, respectively, since they represent the marginal value of real and reactive power at a specific node and hour.
We note that the linearization in \eqref{eq13} is performed around the optimal operating point obtained by the exact AC OPF model Opt1, and that it relates variable $\Delta I_{ij}$ to branch flow model variable $l_{ij}$.
We solve Opt2 to derive dual variables $\lambda _{j}^{P}$ and $\lambda _{j}^{Q}$ (LMVs).
An equivalent approach would be to employ sensitivity analysis, following \cite{CaramanisEtAl2016}, which would require the calculation of the partial derivatives of the branch flow variables \emph{w.r.t.} to the real and reactive power net demand, at the system's optimal operating point.\footnote{
The P-LMV (Q-LMV) at a specific node can be obtained by the partial derivative of the objective function in \eqref{eq12} \emph{w.r.t.} net real (reactive) demand at that node.
For the first two terms we refer to \cite{CaramanisEtAl2016} and the analysis in \cite{Papavasiliou2018}.
The third term involves the partial derivative of $\Delta I_{ij} =\sqrt{l_{ij}} - I_{ij}^{\text{max}}$ which relates to the partial derivative of variable $l_{ij}$ with the coefficient $0.5 \left( \sqrt{l_{ij}^{0}} \right)^{-1}$; see also \eqref{eq13}. 
We also clarify that by measuring the overload in Amps, using variable $\Delta I_{ij}$, we naturally relate the MCC (measured in \$ per Amp) to the upgrade of a line that is typically measured in Amps.
Another option would be to measure the overload in Amps$^2$, and adjust the MCC accordingly.
We would not then need the linearization in \eqref{eq13}, as we could use a variable $\Delta l_{ij} = \max \left[ 0, l_{ij} -{ \left( I_{ij}^{\text{max}} \right)}^2 \right]$. 
This option could be viewed as measuring the overload with the amount of thermal losses above the rated capacity.}

\subsection{Generic DER Procurement} \label{DERprocurement}

In this step, we derive an optimal generic DER allocation that alleviates overload at a specific hour.
We introduce variables $P_{j}^{\text{DER}} \ge 0$, and $Q_{j}^{\text{DER}}$ for real and reactive power procured from generic DERs at node $j$, at a cost equal to P-LMV and Q-LMV, respectively, as estimated in the pricing step.
The new objective function is defined by
\begin{equation}
 \underset{\begin{smallmatrix} P_0, Q_0, P_{ij}, Q_{ij}, \\ 
 v_i, l_{ij}, P_{j}^{\text{DER}}, Q_{j}^{\text{DER}} \end{smallmatrix}}{\mathop{\min}}
 c^P P_0 + c^Q Q_0 + \sum\limits_{ j\in \mathcal{N}^+ }
{\left( \lambda _{j}^{P} P_{j}^{\text{DER}} + \lambda _{j}^{Q} Q_{j}^{\text{DER}} \right)},	\label{eq15}   
\end{equation}
where the time index is omitted since all variables/parameters refer to a specific hour. 
Note that $\lambda _{j}^{P}$ and $\lambda _{j}^{Q}$ are parameters whose values are obtained from the solution of Opt2. 
The power balance constraints \eqref{eq2c}--\eqref{eq2d} are modified accordingly:
\begin{subequations}
\begin{align}
 P_{ij} - l_{ij} r_{ij} + P_j + P_{j}^{\text{DER}} - \sum\limits_{k: j \to k} P_{jk} = 0, \,(\lambda _{j}^{P}) \, \forall j\in \mathcal{N}^{+}, \label{eq16a}	\\
 Q_{ij} - l_{ij} x_{ij} + Q_j + Q_{j}^{\text{DER}}- \sum\limits_{k: j \to k} Q_{jk} = 0, \,(\lambda _{j}^{Q}) \, \forall j\in \mathcal{N}^{+}. \label{eq16b}
\end{align}
\end{subequations}
Network constraints –-- e.g., service transformer rated capacities ---  may impose a bound on the real and reactive power DER quantities that can be procured at a certain node:
\begin{subequations}
\begin{align}
 P_{j}^{\text{DER}}\le \bar{P}_{j}^{\text{DER}}, \,\, \forall j\in \mathcal{N}^{+}, \label{eq17a}	\\
-\bar{Q}_{j}^{\text{DER}}\le Q_{j}^{\text{DER}}\le \bar{Q}_{j}^{\text{DER}}, \,\, \forall j\in \mathcal{N}^{+}. \label{eq17b}
\end{align}
\end{subequations}

The optimal generic DER allocation is obtained by solving the following (QCP/SOCP) optimization problem:
\begin{equation}
\begin{split}
 \textbf{Opt3: } \eqref{eq15}, \text{ s.t. } \eqref{eq2a}-\eqref{eq2b},  \eqref{eq16a}-\eqref{eq16b},  \eqref{eq3}-\eqref{eq6},\\
 \text{ and }  \eqref{eq17a}-\eqref{eq17b}.	\label{eq18}    
\end{split}
\end{equation}
The solution of Opt3 provides an estimate of the DER quantities required to satisfy ampacity constraints at a minimal procurement cost.
In the absence of DER quantity bound constraints \eqref{eq17a}--\eqref{eq17b}, the solution of Opt3 is a lower bound on the actual DER procurement cost. 
Inclusion of constraints \eqref{eq17a}--\eqref{eq17b}, calibrated appropriately for a specific feeder, yields a more realistic estimate of the DER procurement cost.
An advantage of the suggested optimal DER procurement is that all network constraints are observed eliminating the potential of excessive DER additions at one or more locations introducing new problems in back flow, high voltage, etc.

\section{Test Cases} \label{TestCases}

Both test cases were adapted from actual feeders in the ComEd area, IL, representing two typical investment projects. 
For the purposes of this paper, we sanitized the data, while preserving the salient features of the topology and electrical properties, and we employed a high fidelity single-phase AC OPF model.\footnote{
We used the positive sequence of balanced three-phase versions and we compared with three-phase load flow results of the unbalanced feeders.
Since both feeders did not exhibit over/under-voltage issues that might require upgrades targeted to deal with voltage violations --- in which cases potentially high unbalances would require a three-phase representation, the single-phase model proved adequate in illustrating the proposed framework in typical and most representative feeders experiencing overload, in an easy to follow and yet sufficiently realistic and accurate exposition.}
The distribution utility expects load growth and/or potential new customers/loads that absent a DER solution would require a wires investment. 
The cost of this investment can be either associated directly to feeder lines and equipment (Feeder 1) or involve new reconfiguration capability to connect to another feeder (Feeder 2).
We note that both feeders have loop capabilities and tie switches, but they are typically operated in a radial topology through predetermined schemes.
Indeed, network reconfiguration is applied to relieve congestion and mitigate unbalances in the operational timescale.
For the purposes of this paper, topology configuration choices are implicitly captured, since the SOCP model can be applied for different network configurations, allowing for the optimal network topology to be used for each time period, driven by the anticipated loads.\footnote{
An extension of the SOCP problem to explicitly include reconfiguration options, following the formulation proposed in \cite{BaiEtAl2018} and resulting in a Mixed Integer SOCP (MISOCP) problem is straightforward, and it affects only the pre-processing step. 
The MISOCP model, optimizing available reconfiguration actions, can provide the optimal switch settings that yield an SOCP problem reflecting optimal network configuration for a specific load level.
Once the optimal configuration is found, it is passed to the pricing step to calculate LMVs; the pricing step can be applied as is.}

The data for the two feeders are listed next.

\subsection{Feeder 1 (88 Nodes)}
\begin{figure}[tb]
\centering
\includegraphics[width=2.5in]{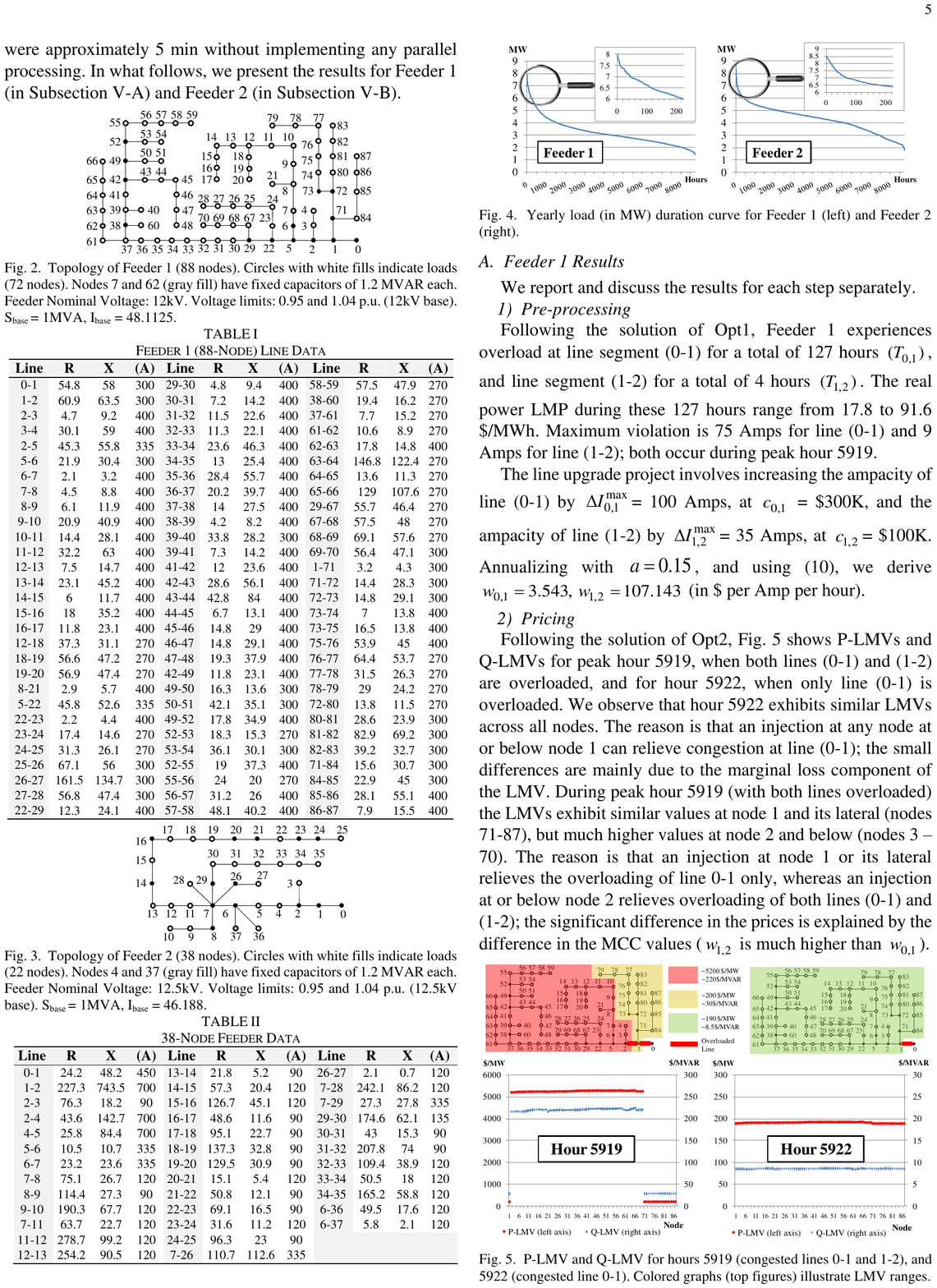}
\caption{Topology of Feeder 1 (88 nodes). Circles with white fills indicate loads (72 nodes). Nodes 7 and 62 (gray fill) have fixed capacitors of 1.2 MVAR each. Feeder Nominal Voltage: 12kV. Voltage limits: 0.95 and 1.04 p.u. (12kV base). Sbase = 1MVA, Ibase = 48.1125.}
\label{Fig2}
\end{figure}

Increased total load is expected to result in line overload close to the root of this 88-node feeder whose topology is shown in Fig. \ref{Fig2}.
Table \ref{tab1} depicts line resistance and reactance (R and X in $10^{-3}$ Ohm), and ampacity (A in Amps).
Line numbering indicates the from- and to-node.
The best investment involves the reinforcement of line segments (0--1) and (1--2) at a total project cost $C = \$400$K.
\begin{table}[tb] 
	\caption{Feeder 1 Line Data}  \label{tab1} 
\centering
\includegraphics[width=3.45in]{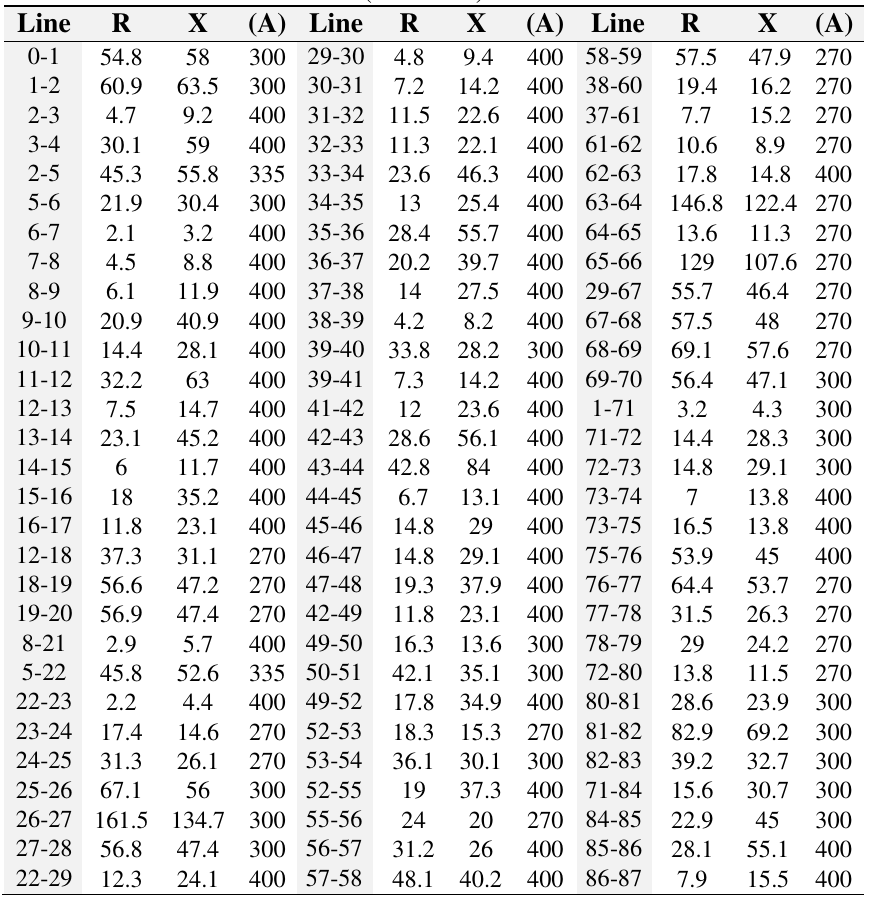}
\end{table}
\begin{figure}[h!]
\centering
\includegraphics[width=2in]{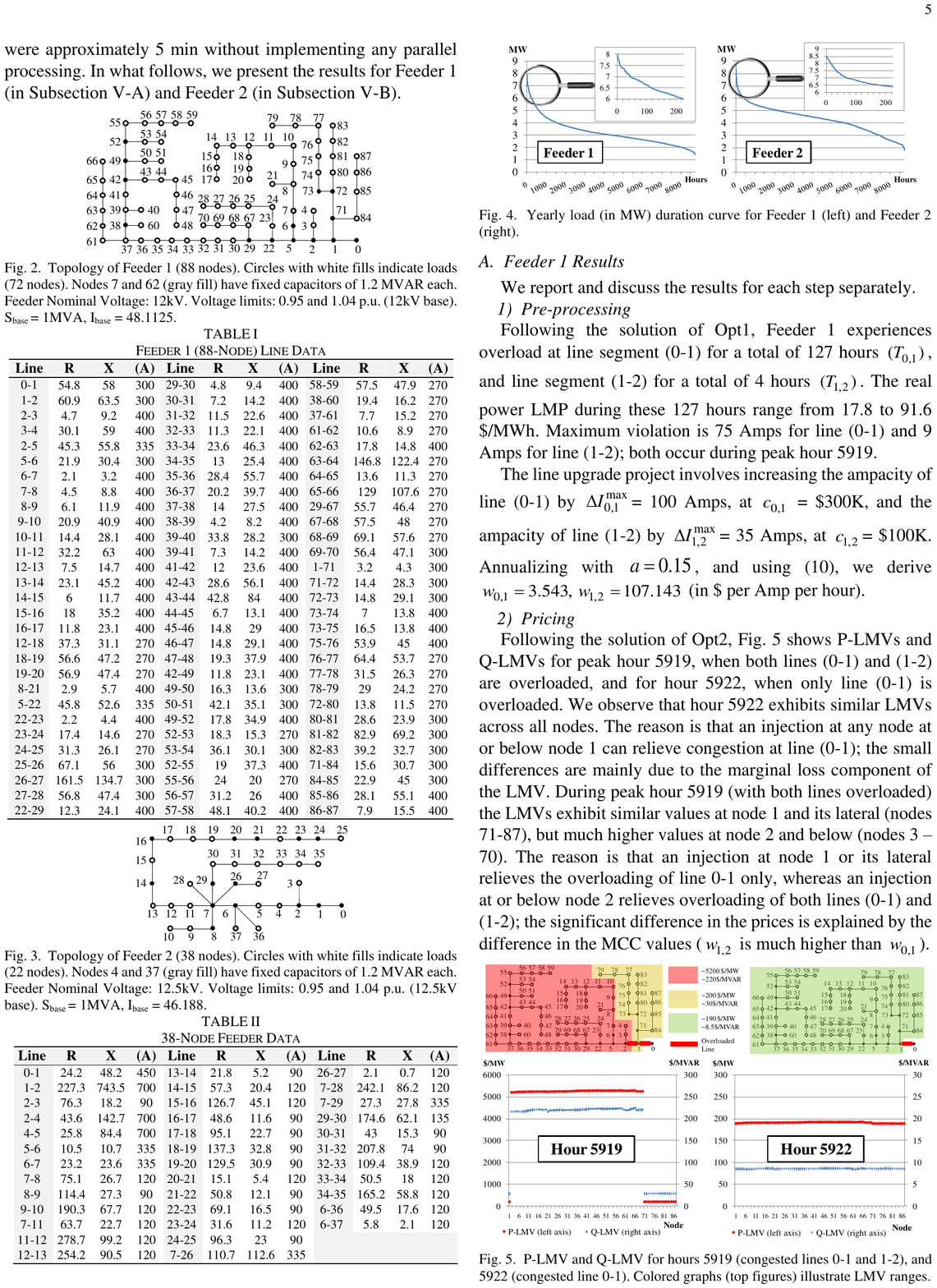}
\caption{Topology of Feeder 2 (38 nodes). Circles with white fills indicate loads (22 nodes). Nodes 4 and 37 (gray fill) have fixed capacitors of 1.2 MVAR each. Feeder Nominal Voltage: 12.5kV. Voltage limits: 0.95 and 1.04 p.u. (12.5kV base). Sbase = 1MVA, Ibase = 46.188.}
\label{Fig3}
\end{figure} 

\subsection{Feeder 2 (38 Nodes)}
Feeder 2 has 38 nodes and is expected to exhibit overload in various lines.
Its topology is shown in Fig. \ref{Fig3} and the line data in Table \ref{tab2}.
The best alternative project for Feeder 2 involves a connection with neighboring feeders, with $C = \$1$M.

\begin{table}[tb] 
	\caption{Feeder 2 Line Data}  \label{tab2} 
\centering
\includegraphics[width=3.45in]{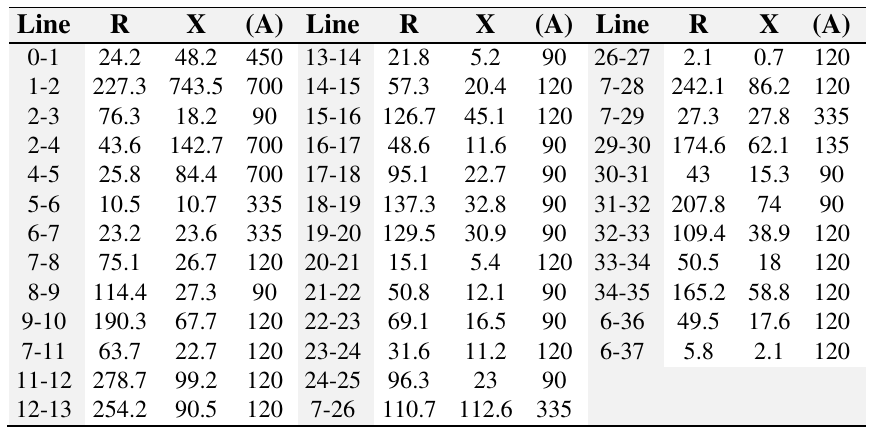}
\end{table}
\begin{figure}[ht]
\centering
\includegraphics[width=3in]{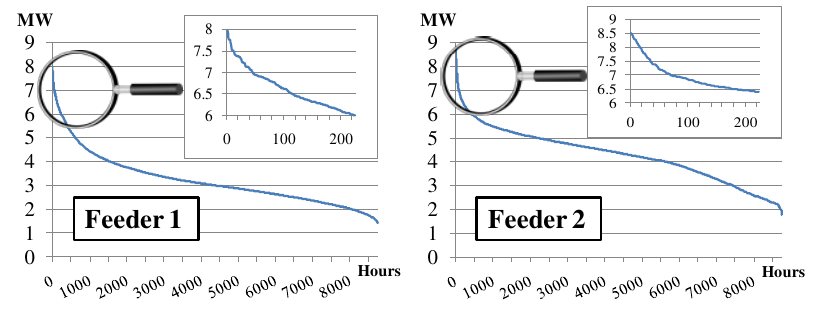}
\caption{Yearly load (in MW) duration curve for Feeder 1 (left) and Feeder 2 (right).}
\label{Fig4}
\end{figure} 

\section{Numerical Results} \label{Results}

Our framework is demonstrated on the two aforementioned test cases.
Yearly load duration curves for both feeders are shown in Fig. \ref{Fig4}.
Power factors at individual nodes range from 0.85 (for commercial nodes) to 0.95 (for residential nodes). 
We use 2016 PJM LMPs in the ComEd area, and assume that the opportunity cost of reactive power is 5\% of the (real power) LMP.
Annualization is done with $\alpha = 0.15$.

The runs were made on a Dell Intel Core i7-5500U @ 2.4 GHz with 8GB RAM, using CPLEX 12.7.
Computational times were approximately 5 min without any parallel processing.
In what follows, we present the results for Feeder 1 (in Subsection \ref{Results1}) and Feeder 2 (in Subsection \ref{Results2}).

\subsection{Feeder 1 Results} \label{Results1}
We report and discuss the results for each step separately.

\subsubsection{Pre-processing}
Following the solution of Opt1, Feeder 1 experiences overload at line segment (0--1) for a total of 127 hours $(T_{0,1})$, and line segment (1--2) for a total of 4 hours $(T_{1,2})$. 
The real power LMP during these 127 hours ranges from 17.8 to 91.6 \$/MWh.
Maximum violation is 75 Amps for line (0--1) and 9 Amps for line (1--2); 
both occur during peak hour 5919.

The line upgrade project involves increasing the ampacity of line (0--1) by $\Delta I_{0,1}^{\max }= 100$ Amps, at $c_{0,1} =\$300$K, and the ampacity of line (1--2) by $\Delta I_{1,2}^{\max }= 35$ Amps, at $c_{1,2} = \$100$K. 
Annualizing and using \eqref{eq10}, we derive $w_{0,1} = 3.543, \, w_{1,2} = 107.143$ (in \$ per Amp per hour).

\subsubsection{Pricing}
\begin{figure}[tb]
\centering
\includegraphics[width=3.45in]{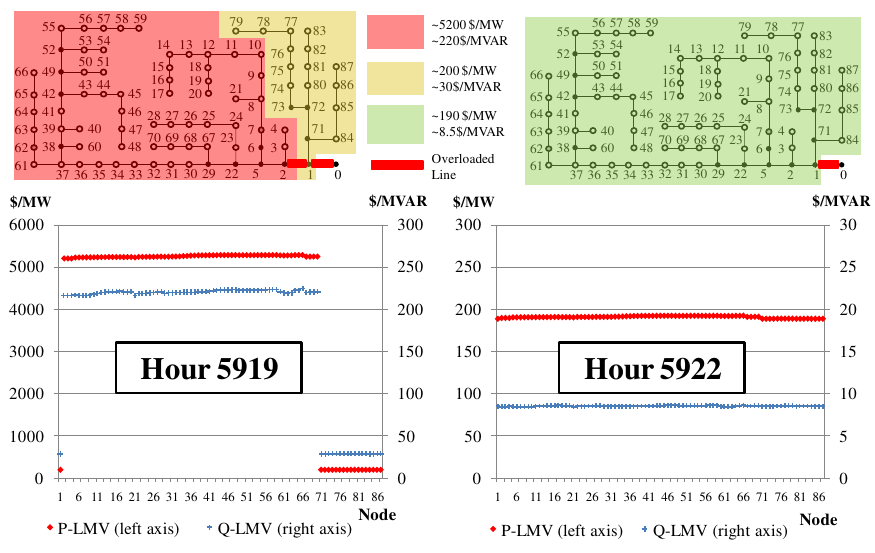}
\caption{P-LMV and Q-LMV for hours 5919 (congested lines 0--1 and 1--2), and 5922 (congested line 0--1). Colored graphs (top figures) illustrate LMV ranges.}
\label{Fig5}
\end{figure}

Following the solution of Opt2, Fig. \ref{Fig5} shows P-LMVs and Q-LMVs for peak hour 5919, when both lines (0--1) and (1--2) are overloaded, and for hour 5922, when only line (0--1) is overloaded.
We observe that hour 5922 exhibits similar LMVs across all nodes.
The reason is that an injection at any node at or below node 1 can relieve congestion at line (0--1);
the small differences are mainly due to the marginal loss component of the LMV.
During peak hour 5919 (with both lines overloaded) the LMVs exhibit similar values at node 1 and its lateral (nodes 71--87), but much higher values at node 2 and below (nodes 3–-70). 
The reason is that an injection at node 1 or its lateral relieves the overload of line 0--1 only, whereas an injection at or below node 2 relieves overload of both lines (0--1) and (1--2);
the significant difference in the prices is explained by the difference in the MCC values ($w_{1,2}$ is much higher than $w_{0,1}$).

\subsubsection{Generic DER Procurement}
In this step, we consider the solution of Opt3 under two scenarios: 
(a) without \eqref{eq17a}--\eqref{eq17b}, i.e., generic DERs are allowed at any node (except the root) and at any capacity, and 
(b) with \eqref{eq17a}--\eqref{eq17b} limiting DERs to the 72 load nodes, and to a maximum of 40 kW/kVAR.

\begin{figure}[tb]
\centering
\includegraphics[width=3.45in]{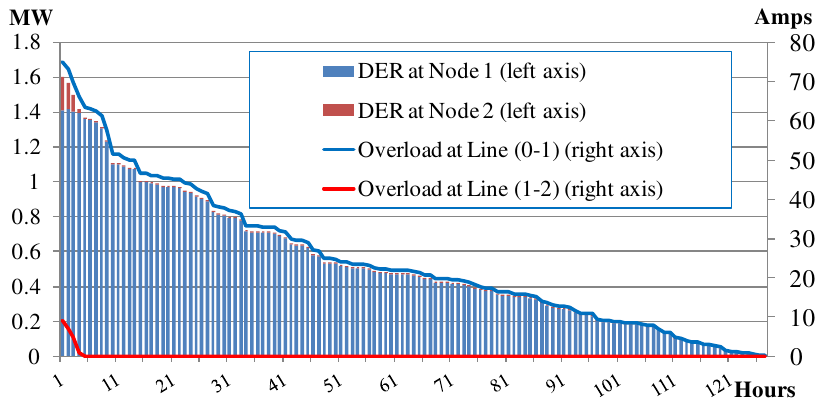}
\caption{Unconstrained scenario. Generic DER real power procurement (nodes 1, 2) for 127 hours of overload (in descending order of overload at line 0--1).}
\label{Fig6}
\end{figure}

We illustrate the results for the unconstrained scenario in Fig. \ref{Fig6}, for the 127 hours of overload.
Fig. \ref{Fig6} shows the DER real power procurement (reactive power is low and its cost negligible) at nodes 1 and 2 (stacked diagram, left axis), and the overload at lines (0--1) and (1--2) in Amps (right axis).
The results confirm that DERs are procured exactly downstream from the overloaded lines.
Note that at the peak hour about 1.6 MW of DER are required to relieve overloads.

\begin{figure}[tb]
\centering
\includegraphics[width=3.45in]{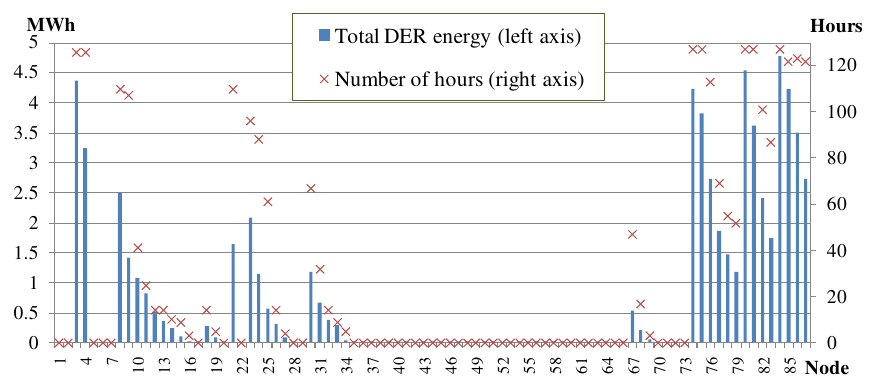}
\caption{Constrained scenario. Generic DER procurement (total energy and hours of energy provision) per node.}
\label{Fig7}
\end{figure}

For the more interesting constrained scenario, we show the generic DER procurement per node in Fig. \ref{Fig7} (total energy and number of hours). 
We observe that most of the DER energy is procured from the first lateral (nodes 74--87), then the second lateral (nodes 3--4), then the third (nodes 8--21), the fourth (nodes 23--28) and the fifth (nodes 67--70).
We also note that DER energy (primarily) and the number of hours (with some exceptions) decrease as we move down a lateral while the LMVs increase as a result of higher marginal losses.

Comparing the two scenarios, the amount of real and reactive power of generic DERs required to relieve overload (as obtained from the optimal solution of Opt3) does not exhibit significant differences: 67.59 MWh and 3.1 MVARh in the unconstrained scenario; 67.29 MWh and 4.6 MVARh in the constrained scenario.
The difference, however, becomes very significant when comparing the generic DER procurement costs. 
For the unconstrained scenario, the procurement cost of real power is \$13,475, and \$38 for reactive power.
For the constrained scenario, the cost increases to \$29,890 and \$84 for real and reactive power, respectively. 
The source of the difference is mainly the 4 hours when both lines ampacities are binding. 
In the constrained scenario, a significant amount of energy (downstream of node 2) is procured at a much higher price (refer to LMVs in Fig. \ref{Fig5}). 
But despite the much higher cost of the constrained scenario, the results indicate that DERs as NWAs are more favorable compared to the wires solution, which would result in a \$60K annual addition to the rate base.

\subsection{Feeder 2 Results} \label{Results2}

\subsubsection{Pre-processing}

Feeder 2 experiences overload during 136 hours.
In particular, 75 hours on line segment (5--6), 75 hours on (29--30), 109 hours on (30--31), 87 hours on (31--32), and 136 hours on (6--36); the maximum overload at the peak hour (5534) is 97.5 Amps (29\% above ampacity), 36.8 (27.3\%), 28 (31.1\%), 25.9 (28.7\%), and 40.3 (33.6\%), respectively. 
The LMPs across the 136 hours range from 17.18 to 72.22 \$/MWh.

Since the investment is part of a reconfiguration project, we allocate the project cost to each line using \eqref{eq11}. 
For line lengths (in m), $L_{5,6} = 100,\, L_{29,30} = 450, \, L_{30,31} = 120, \, L_{31,32} = 550, \, L_{6,36} = 150$, from \eqref{eq10}, we obtain $w_{5,6} = 4.005,\, w_{29,30} = 18.023,\, w_{30,31} = 3.307, \, w_{31,32} = 8.835,\, w_{6,36} = 3.313$.

\subsubsection{Pricing}
\begin{figure}[tb]
\centering
\includegraphics[width=3.45in]{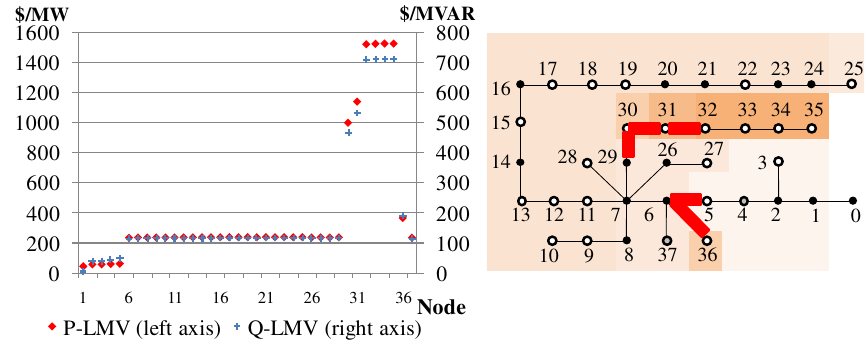}
\caption{P-LMV and Q-LMV for hour 5534; congested lines (5--6), (29--30), (30--31), (31--32), (6--36).}
\label{Fig8}
\end{figure}

Following the solution of Opt2, we illustrate P-LMVs and Q-LMVs for peak hour 5534 in Fig. \ref{Fig8}.
We observe that LMVs increase along the overloaded lines: they exhibit their first step at node 6, and its lateral node 37, then they increase at lateral 36, and take similar values from node 7 to 29, then they increase gradually over nodes 30, 31 and 32, and take similar values at nodes 33-35.
Note that Feeder 2 exhibits lower maximum values of LMVs compared to Feeder 1.
The reason is that Feeder 1 has a much higher MCC value at line (1--2); even though the project cost of Feeder 2 is much higher, there is no line that is overloaded for very few hours --– as is line (1--2) of Feeder 1 --- and hence, the cost is spread across more hours.

\subsubsection{Generic DER Procurement}
\begin{figure}[tb]
\centering
\includegraphics[width=3.45in]{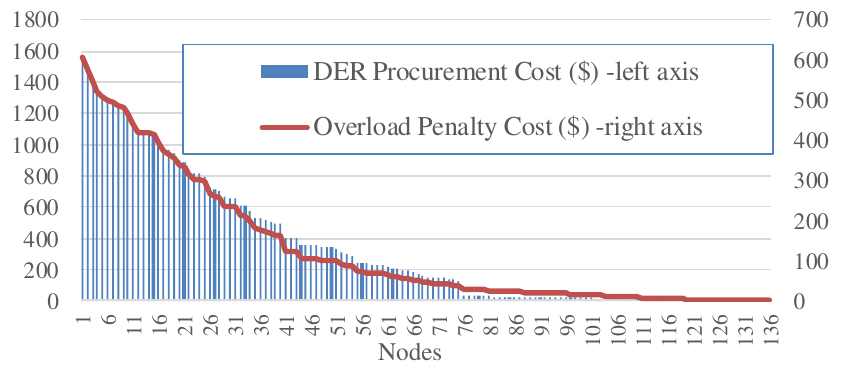}
\caption{Generic DER procurement costs and aggregate overload (constraint violation) cost. Hours are sorted in descending order of DER cost.}
\label{Fig9}
\end{figure}

In Fig. \ref{Fig9}, we plot the generic DER procurement cost for each hour (as obtained from the solution of Opt3), and we see that its pattern generally follows the overload cost, $\sum\nolimits_{(i,j)} w_{ij} \Delta I_{ij}$.
The total procurement cost for real power is \$37,063, and for reactive power \$9,211, for a total of \$46,274. 
Again, compared to the \$150K annual rate base burden of the grid investment, results indicate that DERs constitute an attractive alternative.

\begin{figure}[tb]
\centering
\includegraphics[width=2.7in]{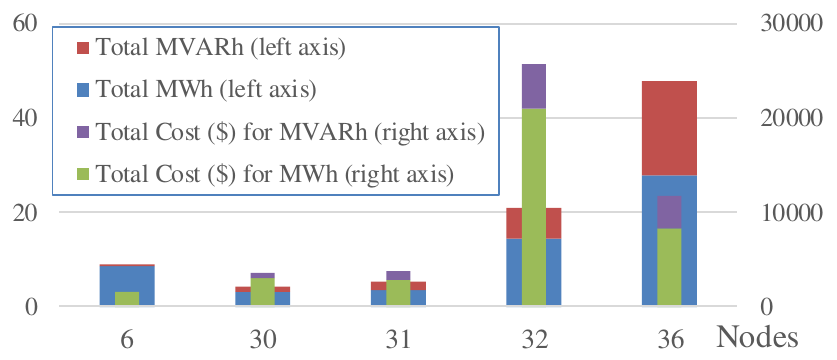}
\caption{Total generic DER procurement (MWh and MVARh) and total generic DER procurement cost per node.}
\label{Fig10}
\end{figure}

However, not surprisingly, location plays a crucial role in this feeder.
More specifically, it is not a surprise that overload of line (6--36) can only be relieved by DERs located at node 36.
On the other hand, overload of line (5--6) can be relieved by DERs located at any node downstream of 6 (similarly to a root node line overload), and overload of line (29--31) can be relieved by nodes 30 (partially), 31 (partially) and downstream of 31 (fully).
Fig. \ref{Fig10} illustrates the total optimal real and reactive power procurement by node along with the procurement costs. 
It is evident that DERs located at node 32 exhibit a much higher (4 to 5 times higher) procurement cost (and hence value) per unit of MWh/MVARh compared to DERs located at node 36. 
Depending on the economic viability of DERs, the results of this analysis could incentivize the planner to consider hybrid solutions, with a partial grid investment, and a partial NWA, that might have not been considered in the absence of this analysis.

\section{Policy Implications and Extensions} \label{Policy}

As already stated, a key objective of the proposed framework has been to enable a distribution utility to rely on information that is in its planning province;
determining the cost of the best grid investment alternative is within the utility's domain and expertise.
For instance, the planner should not have to make assumptions about the costs of new DERs which hypothetically could be procured or incented to be developed in order to relieve constraints and achieve feasibility.

An issue arises that in many cases the best wires alternative may be too ``large'' or too ``lumpy'' to be economic when DER investment alternatives are considered. 
Said differently, if the full cost of a large investment justified by economies of scale and higher future capacity were to be used to value DERs, then the DERs would be overvalued by unjustifiably high overload costs.
Our framework provides two distinct remedies:
First, the cost of the investment is annualized, i.e., the wires investment cost is translated to its annual impact on the rate base.
Second, its cost is pro-rated to the capacity that load growth indicates will be required during the next year or the relevant planning horizon.
Annualization and pro-rating introduces the notion of the MCC, which we use in the valuation of generic DERs that are in fact invariant of actual DER costs and capabilities.

Discussions abound on whether compensation should be locked in, e.g., through grandfathering, to protect early DER adopters from value tanking due to future overexpansion or location specific load growth that renders their DER investments stranded assets.
These discussions are rather focused on the mechanics of DER compensation over time, grandfathering and risk allocation policies;
they have been perennial and inconclusive.
We will not address them as part of this paper, since we do not consider them part of the valuation question that this paper is trying to contribute to.

We wish to comment briefly on one alternative process to the proposed MCC and AC OPF framework which has been on the table for years and relies on the cost of unserved energy.
This process would allow for load reduction via involuntary curtailment at a very high curtailment cost, and would also yield a feasible solution to the DER valuation even when constraints are violated.
It would result, however, in a DER valuation that was very sensitive to the value of lost load and introduce an arbitrary driver of the DER value. 
This would almost certainly become a bone of contention in the planning process and the DER compensation discussion, not to mention that it might be interpreted as a violation of the utility's obligation to serve the load.

The proposed framework has the advantage that it is economically consistent, tied to an objective assessment of grid capacity costs, not related to new DER costs, and not requiring the utility to determine which DERs should be developed, where, and at what cost.
It is not arbitrary and does not result in overstated DER valuations. 
The DER valuation process ought to be focused on the ``value-to-the-grid,'' i.e., the theoretical amount that a utility should compensate a DER owner for its value-to-the-grid, and not solely on revenues from the provision of energy and ancillary services for which existing markets can be resorted to.
To this end, when calculating the DER procurement cost, we do not take into account the energy component of the LMV as represented by the LMP established by the wholesale market at the substation or at the head of each distribution feeder. 
We nevertheless assign the value of marginal losses to the grid LMV.
Hence, we implicitly assume that the utility will compensate the DERs for their energy cost based on the LMP.\footnote{
Indeed, we assume that DERs buy and sell energy at the marginal cost of energy to a typical competitive provider or load serving entity that is usually independent of the DSO or the distribution network utility.
On a typical retail electric service bill, the energy supply rate is on average about half or less than the total kWh rate, hence using the time varying LMP is not an unreasonable approximation.
In summary, the DSO burdens customers for distribution network related costs that enter into the rate base, namely, network asset maintenance, variable and fixed asset costs.
It is precisely the impact of DERs on the DSO rate base that LMVs represent. Customers are billed separately for energy and network (transportation or delivery).
LMV-based non wires alternatives payments are related to the network portion of the bill.}
Said differently, we are interested in the ``value-to-the-grid'' component associated with the non-wires solution, which will be complementing existing market schemes associated with energy and/or ancillary services.

The proposed framework can provide a test for actual DER economic viability of a NWA to inform a contemplated and often expensive effort for DER capacity auctions or other means of encouraging NWAs.
Moreover, averaging of LMVs over more or less homogeneous nodes as we have observed in the reported case studies can provide simpler compensation contracts, should that be desired.
The LMV valuation framework bridges the operating short run and long run investment costs, and is a significant component of the cost-benefit analysis of individual DER investors.
The LMV adds on the LMP the spatiotemporal ``value-to-the-grid'' and along with existing market schemes will have a key role in the investment decision of prospective DER owners.
However, the actual procurement scheme may not, in practice, provide the level of compensation that will entice investors/customers to deploy DERs.
If the cost of DER technology has dropped since the existing DERs were installed, the compensation may be more than adequate.
Otherwise, and particularly when DER investors are worried about annual valuations dropping precipitously in future years, it may be insufficient.
In the test cases presented, the calculated generic DER procurement cost was lower than the cost of the wires solution.
In general, the results will depend on the specific case, the profile of the anticipated overload, the topology and feeder characteristics.
However, LMVs will inform on potential DER procurement methods and help to design them so as to match or even exceed by some measure the wires solution cost, when a decision-maker perceives additional benefits from a DER investment.
But even if the incentives prove to be inadequate for entirely deferring the wires investment, they may though redirect the planning process in partial or lighter feeder upgrades that could have not been considered in the context of conventional investment planning efforts.

\begin{figure}[ht]
\centering
\includegraphics[width=2.7in]{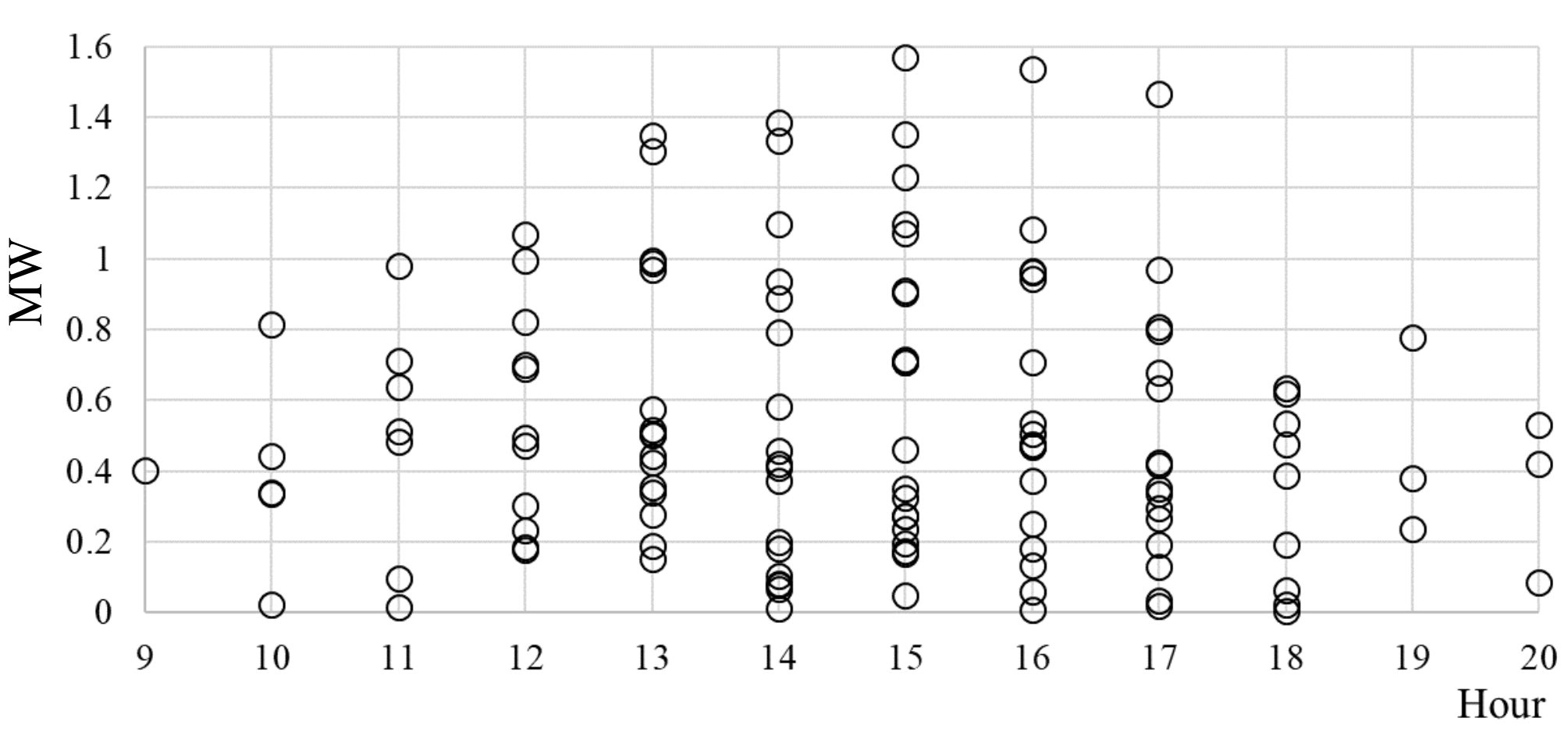}
\caption{Total hourly allocation of real power DER procurement (Feeder 1, constrained scenario).}
\label{Fig11}
\end{figure}

Given the desire to derive as much as possible actual DER-independent NWA results, so far we have stayed away from focusing on actual DERs with their specific capabilities and costs, and we believe that our findings are still useful in introducing the concept of DERs as NWAs.
The P-LMV and Q-LMV of a generic DER at a specific location and hour can be used to calculate the value of an actual DER with specific capabilities.
For instance, a solar PV DER equipped with a smart inverter (assuming it is sized to its nameplate capacity $K$) will be constrained for its real and reactive power provision, $P$ and $Q$, by its capacity, i.e., $P^2 + Q^2 \leq K^2$, and also $P$ will be constrained by the irradiation level (say $\rho$, with $ 0 \leq \rho \leq 1 $), i.e., $P \leq \rho K$. 
The value of this solar PV at each hour will be calculated by the the provided $P$ and $Q$ multiplied with P-LMV and Q-LMV, respectively.
Of course, the hourly allocation of the anticipated overload is significant in determining the ability of a solar PV to act as a NWA, given its irradiation level constraint. 
As an example, we provide in Fig. \ref{Fig11}, the hourly allocation of the overload (in terms of estimated real power required) for the constrained scenario of Feeder 1.
The overload appears in summer daytime hours 9--20, and not surprisingly, solar PV is an excellent fit for contributing in real power as a NWA.
In general, this analysis can be performed for each DER type (even for hybrid systems involving storage), by the utility or the DER investor.

Lastly, while re-conductoring has served as the primary example of wires investments in our case studies, other possibilities such as repowering (raising circuit voltage level), replacing switchgear or limiting station exit cables, and other measures can be similarly treated.
In this respect, the cost of required voltage regulation or circuit impedance reduction, addition of capacitor banks or LTC regulators can be calculated and used to derive appropriate costs for over and under voltage constraint violation.

\section{Conclusions} \label{Conclusions}

We proposed a valuation methodology for DERs as NWAs. 
We employed the concepts of MCC and LMV, and described a framework that uses traditional planning process investment cost information incorporated in an AC OPF problem to derive generic DER LMVs with no need to rely on estimates of actual DER costs and capabilities.
Our framework determines DER values which are locational in space and time, for both real and reactive power.
Generic DER LMVs are invariant to actual DER technology and cost but can be used as the basis for assigning value to and potentially compensate any DER technology, including solar, EVs, demand response, etc.

In our future work, we plan to extend our approach to a three-phase unbalanced system, consider cases of feeders with voltage issues, elaborate on specific grid limitations (e.g., capacity of transformers) and incorporate related costs in the LMVs, and examine the viability of various DER types under certain procurement schemes.
As such, an important issue that remains to be fully investigated includes the reconciliation of ex ante LMV estimates used in our framework to quantify desirable ``generic'' DER procurement with ex post LMVs --- after DERs are in place --- that may lend themselves more appropriately for actual DER compensation.

\section*{Acknowledgment}
The authors gratefully acknowledge the contribution of M. Davoudi in sanitizing the feeders and validating the load flow results, and of F. Farzan for extracting load profiles.

\ifCLASSOPTIONcaptionsoff 
  \newpage
\fi

\end{document}